\documentclass{article}
\usepackage{hyperref}
\usepackage{amsmath,amssymb,amsthm,authblk}
\usepackage[pdftex]{graphicx}

\title{\bf\large Well-posedness and stability for a generalized micropolar thermoelastic body with infinite memory}
\date{\vspace{-5ex}}
\author[1]{\small A. Guesmia}
\author[2]{\small J. E. Mu\~noz Rivera}
\author[3]{\small M. Sep\'ulveda}
\author[4]{\small O. Vera Villagr\a'{a}n}
\affil[1]{\small\it Institut Elie Cartan de Lorraine, UMR 7502,
Universit\'e de Lorraine, 3 Rue Augustin Fresnel, BP 45112, 57073 Metz Cedex 03, France.}
\affil[2]{\small\it Department of Mathematics, Federal University of Rio de Janeiro and National Laboratory for Scientific Computation, Petr\'opolis, RJ, Brasil.}
\affil[3]{\small\it DIM \& CI$^2$MA, Universidad de Concepci\'on, Concepci\'{o}n, Chile.}
\affil[4]{\small\it Departamento de Matem\a'{a}tica, Universidad del B\'{\i}o-B\'{\i}o, Concepci\'{o}n, Chile.}
\affil[1]{aissa.guesmia@univ-lorraine.fr}
\affil[2]{rivera@lncc.br}
\affil[3]{mauricio@ing-mat.udec.cl}
\affil[4]{overa@ubiobio.cl}

\begin{document}
\title{\bf\large Well-posedness and stability of a generalized micropolar thermoelastic body with infinite memory}

\maketitle
\begin{abstract}
We study in this paper the well-posedness and stability of a linear system of a thermoelastic Cosserat medium with infinite memory, where the Cosserat medium is a continuum in which each point has the degrees of freedom of a rigid body.
\end{abstract}

{\bf keywords:} Micropolar thermoelasticity; second sound; semigroups; stability; lyapunov functional.

\numberwithin{equation}{section}
\newtheorem{theorem}{Theorem}[section]
\newtheorem{lemma}[theorem]{Lemma}
\newtheorem{remark}[theorem]{Remarks}
\newtheorem{definition}[theorem]{Definition}
\newtheorem{proposition}[theorem]{Proposition}
\newtheorem{cor}[theorem]{Corollary}
\allowdisplaybreaks

\section{Introduction}

In Ferreira {\it et al.} \cite{ferreira}, the authors considered an isotropic generalized thermoelastic body 
$\Omega \subset \mathbb{R}^{3}$ of density $\rho > 0$. In orthogonal curvilinear coordinates at point $(x_{1},\,x_{2},\,x_{3}) = \boldsymbol{x} \in \Omega$ and time $t > 0$, the displacement field 
\begin{equation*}
\boldsymbol{u} (\boldsymbol{x},\,t) = (u_{1}(\boldsymbol{x},\,t),\,u_{2}(\boldsymbol{x},\,t),\,u_{3}(\boldsymbol{x},\,t) )\in \mathbb{R}^{3},
\end{equation*}
the microrotational field 
\begin{equation*}
\boldsymbol{\varphi}(\boldsymbol{x},\,t) = (\varphi_{1}(\boldsymbol{x},\,t),\,\varphi_{2}(\boldsymbol{x},\,t),\,\varphi_3 (\boldsymbol{x},\,t) )\in \mathbb{R}^{3}
\end{equation*}
and the relative absolute temperature $\theta(\boldsymbol{x},\,t) \in \mathbb{R}$ satisfy the following system (for the generalized thermoelasticity, see, for example, \cite{iesan}, chap. 5, p. 165):
\begin{align*}
\begin{cases}
& \rho\,\boldsymbol{u}_{tt} - (\mu + \kappa)\,{\rm{div}}\,\nabla \boldsymbol{u} - (\lambda + \mu)\,\nabla\,{\rm{div}}\,\boldsymbol{u} - \kappa\,(\nabla\times\boldsymbol{\varphi}) + b\,\nabla \theta + b\,\alpha_{0}\,\nabla\theta_{t} = \rho\,{\boldsymbol{f}} + \boldsymbol{F}, \\
& \rho\,j\,\boldsymbol{\varphi}_{tt} + 2\,\kappa\,\boldsymbol{\varphi} - \gamma\,{\rm{div}}\nabla \boldsymbol{\varphi} - (\alpha + \beta)\,\nabla\,{\rm{div}}\,\boldsymbol{\varphi} - \kappa\,(\nabla\times \boldsymbol{u}) = \rho\,\boldsymbol{g} + \boldsymbol{G},  \\
& h\,\theta_{tt} + d\,\theta_{t} - \frac{1}{T_{0}}\, {\rm{div}} \left( k\,\nabla \theta \right) + b\,{\rm{div}}\, \boldsymbol{u}_{t} = \frac{\rho}{T_{0}} \, S,
\end{cases}
\end{align*}
where we used bold letter for vector fields and the subscript $t$ for the derivative with respect to $t$. Symbols $\nabla$, $\rm{div}$ and $\nabla \times$ stand, respectively, for the spatial gradient, divergence and curl  operators.
On the right hand side, $\boldsymbol{f}$, $\boldsymbol{g}$ and $S$ represent, respectively, the body force, the body couple and the heat source per unit mass, while $\boldsymbol{F}$ and $\boldsymbol{G}$ are additional forces applied to the system, possibly adding dissipativity that will be made explicitly below (\,see \eqref{def-F_G}\,).
We assume  $j >0$\footnote{This is for a spin-isotropic or microisotropic material in the terminology of Eringen.}. Moreover, the constant $b >0$ represents the coupling to the thermal effect, and $h > 0$ and $\alpha_{0} > 0$ come from the model of heat conduction (with second sound). Note that the case $h = \alpha_{0} = 0$ reduces the heat conduction to the classical Fourier heat effect. Finally, the set of parameters $\mu,\,\kappa,\,\lambda,\,\gamma,\,\alpha,\,\beta,\,d,\,k$ and $T_{0}$ are called the constitutive moduli of the body and satisfy the following compatibility conditions (see, for example, \cite{iesan}, p. 113):
\begin{align}
\label{moduli}
\begin{cases}
& 3\,\lambda + 2\,\mu + \kappa > 0, \quad 2\,\mu + \kappa > 0,\quad \kappa > 0,\\
& 3\,\alpha + \beta + \gamma > 0, \quad \gamma + \beta  > 0, \quad \gamma - \beta > 0,\\
& d >0, \quad  k > 0,\quad T_0 >0.
\end{cases}
\end{align}
The above conditions \eqref{moduli} are required in order, for the so-called internal (potential) energy per unit of volume $E_{w}$ (\,see \eqref{214}-\eqref{216} below for our notations\,), to be a definite positive quadratic form. For more details on these conditions, we refer to \cite{neff1, neff2}. 

\hfill

In the sequel, for sake of simplicity, we consider an homogeneous material, so that all the above parameters are constants.   Moreover, we shall set
\begin{equation}\label{simplification}
\rho = \, j \, = \, h\, = \, T_{0} \, =\, 1,
\end{equation}
and the body forces $\boldsymbol{f}$ and $\boldsymbol{g}$ and the heat source $S$ to zero. Finally we will consider dissipative source terms of frictional and memory types, namely
\begin{equation}
\label{def-F_G}
\boldsymbol{F} = -\ \xi_{0}\,\boldsymbol{u}_{t} \quad\hbox{and}\quad \boldsymbol{G} = -\ \xi\,\boldsymbol{\varphi}_{t} - \int_{0}^{+\infty}g(s)\,\Delta\boldsymbol{\varphi}(\boldsymbol{x},\,t - s)\ ds,
\end{equation}
where $\xi$ and $\xi_{0}$ are positive constants and $g: \mathbb{R}^{+}\rightarrow  \mathbb{R}^{+}$ is a given function satisfying some hypotheses (\,see ({\bf G1}) and ({\bf G2}\,) below). To our knowledge, it is the first time that such a dissipative effect is investigated for this kind of material. Another dissipative effect (viscous term of Kelvin-Voigt type) has been studied though \cite{elkaramani1}. 

\hfill

In this paper, we are then concerned with the following initial boundary value problem with frictional dampings and an infinit memory:
\begin{equation}
\label{101} \boldsymbol{u}_{tt} - (\mu + \kappa)\,\Delta \boldsymbol{u} - (\lambda + \mu)\,\nabla\,{\rm{div}}\,\boldsymbol{u} - \kappa\,(\nabla\times \boldsymbol{\varphi} ) + b\,\big( \nabla\theta +\alpha_{0}\,\nabla\theta_{t} \big) +\, \xi_{0}\,\boldsymbol{u}_{t} =0,  
\end{equation}
\begin{equation}
\begin{aligned}
\label{102}& \boldsymbol{\varphi}_{tt} + 2\,\kappa\,\boldsymbol{\varphi} -\gamma\,\Delta\boldsymbol{\varphi} - (\alpha + \beta)\,\nabla\,{\rm{div}}\,\boldsymbol{\varphi} - \kappa\,(\nabla\times \boldsymbol{u}) \\
&+ \int_{0}^{+\infty}g(s)\ \Delta\boldsymbol{\varphi}(\boldsymbol{x},\,t - s)\ ds + \xi\,\boldsymbol{\varphi}_{t} =0,  
\end{aligned}
\end{equation}
\begin{equation}
\label{103} \theta_{tt} + d\,\theta_{t} - k\,\Delta\theta + b\,{\rm{div}}\,\boldsymbol{u}_{t} = 0,
\end{equation}
for any $(\boldsymbol{x},\,t)\in \Omega\times \mathbb{R}^{+}$, where the domain $\Omega$ is assumed to be a simply connected smooth subset of $\mathbb{R}^{3}$ with smooth boundary $\partial\Omega := \Gamma$. The set of equations \eqref{101}-\eqref{103} is completed with homogeneous Dirichlet boundary conditions
\begin{equation}
\label{108}\boldsymbol{u} (\boldsymbol{x},\,t)=\boldsymbol{\varphi} (\boldsymbol{x},\,t)=\theta (\boldsymbol{x},\,t)=0, \quad  \forall\; (\boldsymbol{x},\,t) \in \Gamma \times \mathbb{R}^{+}
\end{equation}
and initial data
\begin{equation}
\begin{aligned}
\begin{cases}
\label{109}& \boldsymbol{u}(\boldsymbol{x},\,0) = \boldsymbol{u}_{0}(\boldsymbol{x}),\quad \boldsymbol{u}_{t}(\boldsymbol{x},\,0) = \boldsymbol{u}_{1}(\boldsymbol{x}),  \\
& \boldsymbol{\varphi}(\boldsymbol{x},\,-\,t) = \boldsymbol{\varphi}_{0}(\boldsymbol{x},\,t),\quad  \boldsymbol{\varphi}_{t}(\boldsymbol{x},\,0) = \boldsymbol{\varphi}_{1}(\boldsymbol{x}),\\
& \theta(\boldsymbol{x},\,0) = \theta_{0}(\boldsymbol{x}),\quad \theta_{t}(\boldsymbol{x},\,0) = \theta_{1}(\boldsymbol{x}),
\end{cases}
\end{aligned}
\end{equation}
where $(\boldsymbol{x},\,t)\in \Omega\times \mathbb{R}^{+}$. The boundary conditions correspond to a rigidly clamped structure with temperature held constant at the extremities (\,equal to $T_{0}$\,). The particular uncoupled equation \eqref{102} (\,Lam\'e system\,) with infinite memory (\,\eqref{102} with $\kappa =\xi =0$\,) was considered in \cite{gubc}. 

\hfill

The steel is the commonest engineering structural material, hence the linear theory of elasticity is of main  importance in the stress analysis of steel. Somehow linear elasticity describes the mechanical behavior of other common solid material, e.g. concrete, wood and coal. However, the theory does not apply to the behavior of many of the new synthetic materials of the elastomer and polymer type, e.g. polymethyl-methacrylate (perspex), polyethylene and polyvinyl chloride. To represent the behavior of such materials, the linear theory of micropolar elasticity is adequate because it takes into consideration the granular character of the medium and it is intended to be applied to materials for which the ordinary classical theory of  elasticity fails owing to the microstructure of the material (see \cite{mindlin}). Within such theory, solids can undergo macro-deformations and micro-rotations. The motion in this kind of solids is completely characterized by the displacement vector and the microrotational vector, whereas in case of classical elasticity, the motion is characterized by the displacement vector only. Basically, the difference between classical continuum theories and the micropolar theory is that the latter admits independent rotations of the material's substructure; that is the local intrinsic rotations (\,microrotations\,) which are taken to be kinematically independent of the linear displacements. It is believed that such theory is applicable in the treatment of granular and fibrous composite materials. 

\hfill

The micropolar theory have been extended to include thermal effects by Eringen \cite{eringen, eringen1} and Nowacki \cite{nowacki2, nowacki1, nowacki3, nowaki}. A generalized theory of micropolar thermoelasticity was extended in \cite{I}. In the last years, the theory of thermoelasticity for bodies with microstructure has been deserved much attention. A thermodynamic theory for elastic materials with inner structure whose particles, in addition to microdeformations, possess microtemperatures was proposed in \cite{grot}. A theory of micromorphic fluids was developed by \cite{riha1, riha2}. Different types of problems in micropolar thermoelasticity have been studied by several authors, see, for example, \cite{eringen1, eringen3, eringen2, iesan0}. The linear theory of thermoelasticity with microtemperatures for materials with inner structure whose particles, in addition to the classical displacement and temperature fields, possess microtemperaturas was presented in \cite{iq}, where an existence theorem was proved and the continuous dependence of solutions of the initial data and body loads was established. In this article, we consider field equations in terms of the displacement vector, microrotation vector and temperature variation. We investigate in sections 2 and 3 the well-posedness of \eqref{101}-\eqref{109} and the asymptotic behavior of its solutions as $t\rightarrow +\infty$. We finish our paper by giving some general comments and issues in section 4.

\section{Assumptions and setting of the semigroup}

In this section, we consider some hypotheses on the relaxation function $g$ and establish the well-posedness of system
\eqref{101}-\eqref{109}. We denote by $L^{2}(\Omega)$ the classical space of square integrable functions over $\Omega$ with an abuse of writing for vector valued functions. Let $\Vert \cdot \Vert$ be the standard $L^2$ norm over $\Omega$ generated by the classical inner product 
\begin{equation*}
\left\langle u,\,{\tilde {u}}\right\rangle_{L^{2}(\Omega)}= \int_{\Omega}u(\boldsymbol{x})\ {\tilde {u}}(\boldsymbol{x})\ d\boldsymbol{x}.
\end{equation*}
Recall also $H^{1}_{0}(\Omega)$, the classical homogeneous Hilbert space endowed with its inner product
\begin{equation*}
\left\langle u ,{\tilde {u}}\right\rangle_{H^{1}_{0} (\Omega)} =  \int_{\Omega} \,\nabla u(\boldsymbol{x})\cdot\nabla {\tilde {u} (\boldsymbol{x})}\ d\boldsymbol{x}.
\end{equation*}
To guarantee the well-posednes of system \eqref{101}-\eqref{109}, we consider the following hypotheses:
\\ 
\\
{\bf (G1)} We assume that $g: \mathbb{R}^{+}\rightarrow  \mathbb{R}^{+}$ is a nonincreasing differentiable function such that
\begin{equation}
\label{201} \int_{0}^{+\infty}g(s)\ ds := \ g_{0} < \gamma 
\end{equation}
and there exists $d_{0}>0$ such that
\begin{equation}
\label{201*} g'(s)\geq -\,d_{0}\,g(s),\quad \forall\; s\in \mathbb{R}^{+}. 
\end{equation}
Moreover, we assume that \eqref{moduli} holds such that
\begin{equation}
\label{hyp}
\alpha_{0}\,d - 1 > 0.
\end{equation}
We thus define the phase space associated with our set of equations \eqref{101}-\eqref{109} by
\begin{equation*}
\mathcal{H}=[H_{0}^{1}(\Omega)]^{3} \times [L^{2}(\Omega)]^{3} \times [H_{0}^{1}(\Omega)]^{3}
\times [L^{2}(\Omega)]^{3}\times H_{0}^{1}(\Omega) \times L^{2}(\Omega)\times [L_g ]^3 ,
\end{equation*}
where $L_{g} = L_{g}^{2} (\mathbb{R}^{+};\, H_{0}^{1}(\Omega))$ is the Hilbert space of all $H_{0}^{1}(\Omega)$-valued and square integrable function defined on the measure space $(\mathbb{R}^{+} \to H_{0}^{1}(\Omega);\, g\ ds)$; that is
\begin{equation}
\label{301}L_{g} = \left\{\eta:\,\mathbb{R}^{+}\rightarrow H_{0}^{1}(\Omega),\quad \int_{\Omega}\int_{0}^{+\infty}g(s)\,|\nabla\boldsymbol{\eta}(\boldsymbol{x},\,s)|^{2}\ ds\ d\boldsymbol{x}<+\infty\right\},
\end{equation}
equipped with the norm that generated by the inner oroduct
\begin{equation*}
\left\langle\eta,\,{\tilde {\eta}}\right\rangle_{L_{g}}= \int_{\Omega}\int_{0}^{+\infty}g(s)\,\nabla \boldsymbol{\eta}(\boldsymbol{x},\,s)\cdot\nabla {\tilde {\boldsymbol{\eta}}}(\boldsymbol{x},\,s)\ ds\ d\boldsymbol{x}.
\end{equation*}
In order to rewrite equation \eqref{101}, we note that 
\begin{equation}
\label{101*}\Delta \boldsymbol{u} = \nabla{\rm{div}}\,\boldsymbol{u} - \nabla \times \nabla \times \boldsymbol{u}.
\end{equation}
Indeed, for $\boldsymbol{u}(x,\,y,\,z)= (u_{1}(x,\,y,\,z),\,u_{2} (x,\,y,\,z),\,u_{3}(x,\,y,\,z))$, we have
\begin{equation*}
\nabla\,div\,u = \nabla\left(\frac{\partial u_{1}}{\partial x} + \frac{\partial u_{2}}{\partial y} + \frac{\partial u_{3}}{\partial z}\right)  
\end{equation*}
\begin{equation*}
=\left(\frac{\partial}{\partial x}\left(\frac{\partial u_{1}}{\partial x} + \frac{\partial u_{2}}{\partial y} + \frac{\partial u_{3}}{\partial z}\right),\,\frac{\partial}{\partial y}\left(\frac{\partial u_{1}}{\partial x} + \frac{\partial u_{2}}{\partial y} + \frac{\partial u_{3}}{\partial z}\right),\,\frac{\partial}{\partial z}\left(\frac{\partial u_{1}}{\partial x} + \frac{\partial u_{2}}{\partial y} + \frac{\partial u_{3}}{\partial z}\right)\right)  
\end{equation*}
\begin{equation}
\label{102*}=\left(\frac{\partial^{2}u_{1}}{\partial x^{2}} + \frac{\partial^{2}u_{2}}{\partial x\partial y} + \frac{\partial^{2}u_{3}}{\partial x\partial z},\,\frac{\partial^{2}u_{1}}{\partial y\partial x} + \frac{\partial^{2}u_{2}}{\partial y^{2}} + \frac{\partial^{2}u_{3}}{\partial y\partial z},\,\frac{\partial^{2}u_{1}}{\partial z\partial x} + \frac{\partial^{2}u_{2}}{\partial z\partial y} + \frac{\partial^{2}u_{3}}{\partial z^{2}}\right)
\end{equation}
and
\begin{equation*}
\Delta \boldsymbol{u} = \left(\Delta u_{1},\,\Delta u_{2},\,\Delta u_{3}\right)  
\end{equation*}
\begin{equation}
\label{103*} = \left(\frac{\partial^{2}u_{1}}{\partial x^{2}} + \frac{\partial^{2}u_{1}}{\partial y^{2}} + \frac{\partial^{2}u_{1}}{\partial z^{2}},\,\frac{\partial^{2}u_{2}}{\partial x^{2}} + \frac{\partial^{2}u_{2}}{\partial y^{2}} + \frac{\partial^{2}u_{2}}{\partial z^{2}},\,\frac{\partial^{2}u_{3}}{\partial x^{2}} + \frac{\partial^{2}u_{3}}{\partial y^{2}} + \frac{\partial^{2}u_{3}}{\partial z^{2}}\right).
\end{equation}
Then \eqref{102*} and \eqref{103*} imply that
\begin{equation*}
\nabla\,div\,\boldsymbol{u} - \Delta \boldsymbol{u}
= \Bigl(\frac{\partial^{2} u_{2} }{\partial x\partial y} + \frac{\partial^{2}u_{3}}{\partial x\partial z} - \frac{\partial^{2}u_{1}}{\partial y^{2}} - \frac{\partial^{2}u_{1}}{\partial z^{2}},\,\frac{\partial^{2}u_{1}}{\partial y\partial x} + \frac{\partial^{2}u_{3}}{\partial y\partial z} - \frac{\partial^{2}u_{2}}{\partial x^{2}} - \frac{\partial^{2}u_{2}}{\partial z^{2}},
\end{equation*}
\begin{equation}
\label{104*} \quad \frac{\partial^{2}u_{1}}{\partial z\partial x} + \frac{\partial^{2}u_{2}}{\partial z\partial y} - \frac{\partial^{2}u_{3}}{\partial x^{2}} - \frac{\partial^{2}u_{3}}{\partial y^{2}}\Bigr).
\end{equation}
On the other hand, we see that
\begin{equation*}
\nabla\times \boldsymbol{u} = i\left(\frac{\partial u_{3}}{\partial y} - \frac{\partial u_{2}}{\partial z}\right) + j\left(\frac{\partial u_{1}}{\partial z} - \frac{\partial u_{3}}{\partial x}\right) + k\left(\frac{\partial u_{2}}{\partial x} - \frac{\partial u_{1}}{\partial y}\right).
\end{equation*}
Then
\begin{equation*}
\nabla\times\nabla\times \boldsymbol{u} = \left |\begin{array}{ccc}  
 i & j & k \\
\frac{\partial}{\partial x} & \frac{\partial}{\partial y} & \frac{\partial}{\partial z} \\
\left(\frac{\partial u_{3}}{\partial y} - \frac{\partial u_{2}}{\partial z}\right) & \left(\frac{\partial u_{1}}{\partial z} - \frac{\partial u_{3}}{\partial x}\right) & \left(\frac{\partial u_{2}}{\partial x} - \frac{\partial u_{1}}{\partial y}\right) 
\end{array}\right| 
\end{equation*}
\begin{equation*}
= i\left[\frac{\partial}{\partial y}\left(\frac{\partial u_{2}}{\partial x} - \frac{\partial u_{1}}{\partial y}\right) - \frac{\partial}{\partial z}\left(\frac{\partial u_{1}}{\partial z} - \frac{\partial u_{3}}{\partial x}\right)\right] 
\end{equation*}
\begin{equation*} 
+\ j\left[\frac{\partial}{\partial z}\left(\frac{\partial u_{3}}{\partial y} - \frac{\partial u_{2}}{\partial z}\right) - \frac{\partial}{\partial x}\left(\frac{\partial u_{2}}{\partial x} - \frac{\partial u_{1}}{\partial y}\right)\right] 
\end{equation*}
\begin{equation*}
+\ k\left[\frac{\partial}{\partial x}\left(\frac{\partial u_{1}}{\partial z} - \frac{\partial u_{3}}{\partial x}\right) - \frac{\partial}{\partial y}\left(\frac{\partial u_{3}}{\partial y} - \frac{\partial u_{2}}{\partial z}\right)\right] 
\end{equation*}
\begin{equation*} 
= \Bigl(\frac{\partial^{2}u_{2}}{\partial y\partial x} + \frac{\partial^{2}u_{3}}{\partial z\partial x} - \frac{\partial^{2} u_{1}}{\partial y^{2}} - \frac{\partial^{2}u_{1}}{\partial z^{2}},\,\frac{\partial^{2}u_{1}}{\partial x\partial y} + \frac{\partial^{2}u_{3}}{\partial z\partial y} - \frac{\partial^{2}u_{2}}{\partial x^{2}} - \frac{\partial^{2}u_{2}}{\partial z^{2}},
\end{equation*}
\begin{equation}
\label{105*} \frac{\partial^{2}u_{1}}{\partial x\partial z} + \frac{\partial^{2}u_{2}}{\partial y\partial z} - \frac{\partial^{2}u_{3}}{\partial x^{2}} - \frac{\partial^{2}u_{3}}{\partial y^{2}}\Bigr).
\end{equation}
From \eqref{104*} and \eqref{105*}, \eqref{101*} follows. So, thanks to \eqref{101*}, we can rewrite equation \eqref{101} as follows:
\begin{equation}
\label{104}
\boldsymbol{u}_{tt} - \mu \,\Delta \boldsymbol{u}  + \kappa\, \nabla \times \nabla \times \boldsymbol{u} - (\lambda + \mu + \kappa)\,\nabla\,{\rm{div}}\,\boldsymbol{u} - \kappa\,(\nabla\times \boldsymbol{\varphi} ) + b\,\nabla \big( \theta + \alpha_{0}\,\theta_{t} \big) + \xi_{0}\,\boldsymbol{u}_{t} =0.
\end{equation}
Now, we rewrite equation \eqref{102} in another form. To do so, following the idea given by \cite{dafermos}, we set 
\begin{equation}\label{2050}
\boldsymbol{\eta}(\boldsymbol{x},\,t,\,s) = \boldsymbol{\varphi} (\boldsymbol{x},\,t) - \boldsymbol{\varphi} (\boldsymbol{x},\,t - s),\quad \forall\;\boldsymbol{x}\in \Omega,\ \forall\;t,\,s\in \mathbb{R}^{+} 
\end{equation}
and
\begin{equation}\label{2051}
\boldsymbol{\eta}_0 (\boldsymbol{x},\,s) =\boldsymbol{\varphi}_{0} (\boldsymbol{x},\,0) - \boldsymbol{\varphi}_{0} (\boldsymbol{x},\,s),\quad \forall\;\boldsymbol{x}\in \Omega,\ \forall\;s\in \mathbb{R}^{+}.
\end{equation}
Then
\begin{equation}\label{205}
\begin{aligned}
\begin{cases}
\boldsymbol{\eta}_{t} (\boldsymbol{x},\,t,\,s) + \boldsymbol{\eta}_{s} (\boldsymbol{x},\,t,\,s)- \boldsymbol{\varphi}_{t} (\boldsymbol{x},\,t) = 0, & \quad \forall\;\boldsymbol{x}\in \Omega,\ \forall\;t,\,s\in \mathbb{R}^{+}, \\
\boldsymbol{\eta} (\boldsymbol{x},\,t,\,s)=0 ,&\quad \forall\; \boldsymbol{x}\in \Gamma,\ \forall\;t,\,s\in \mathbb{R}^{+}, \\
\boldsymbol{\eta} (\boldsymbol{x},\,t,\,0)=0 ,&\quad\forall \;\boldsymbol{x}\in \Omega,\ \forall\;t\in \mathbb{R}^{+}, \\
\boldsymbol{\eta} (\boldsymbol{x},\,0,\,s) =\boldsymbol{\eta}_0 (\boldsymbol{x},\,s),&\quad\forall\;\boldsymbol{x}\in \Omega,\ \forall\;s\in \mathbb{R}^{+}.
\end{cases}
\end{aligned}
\end{equation}
Therefore, the equation \eqref{102} can be written as
\begin{eqnarray}
&  & \boldsymbol{\varphi}_{tt} - \beta_{0}\,\Delta\boldsymbol{\varphi} + 2\,\kappa\,\boldsymbol{\varphi} - (\alpha + \beta)\,\nabla\,{\rm{div}}\,\boldsymbol{\varphi} - \kappa\,(\nabla\times \boldsymbol{u})  \nonumber \\
\label{207} &  & - \int_{0}^{+\infty}g(s)\,\Delta\boldsymbol{\eta} (x,\,t,\, s)\, ds + \xi\,\boldsymbol{\varphi}_{t} = 0,
\end{eqnarray}
where $\beta_0 =\gamma -g_0$. So the system \eqref{101}-\eqref{109} is equivalent to \eqref{103}-\eqref{109}, \eqref{104}, \eqref{2050} and \eqref{207}. Let 
\begin{equation*}
\begin{aligned}
\begin{cases}
\boldsymbol{v}=\boldsymbol{u}_{t} ,\quad \boldsymbol{\psi}=\boldsymbol{\varphi}_t ,\quad\Theta =\theta_{t} ,\\ 
U = (\boldsymbol{u},\,\boldsymbol{v},\,\boldsymbol{\varphi},\,\boldsymbol{\psi},\,\theta,\,\Theta,\,\boldsymbol{\eta})^{T},
\\
U_{0} =(\boldsymbol{u}_{0},\,\boldsymbol{u}_{1},\,\boldsymbol{\varphi}_{0},\,\boldsymbol{\varphi}_{1},\,\theta_{0},\,\theta_{1} ,\,\boldsymbol{\eta}_{0} )^{T},
\end{cases}
\end{aligned}
\end{equation*}
where $T$ means the transpose. The set of equations \eqref{103}-\eqref{109}, \eqref{104}, \eqref{205} and \eqref{207} can be written under the form of an abstract first order evolution problem
\begin{equation}\label{302} 
\begin{aligned}
\begin{cases}
U_{t} (\boldsymbol{x},\,t)= \mathcal{A}U (\boldsymbol{x},\,t), 
& \quad \forall\;\boldsymbol{x}\in \Omega ,\ \forall \;t\in \mathbb{R}^{+},  \\
U (\boldsymbol{x},\,0)= U_{0} (\boldsymbol{x}), 
& \quad \forall\;\boldsymbol{x}\in \Omega,
\end{cases}
\end{aligned}
\end{equation}
where $\mathcal{A}$ is an unbounded linear operator $\mathcal{A}:\mathcal{D}(\mathcal{A})\subset{\cal H} \rightarrow{\cal H}$ defined by 
\begin{equation*}
\mathcal{A}U =  \left(\begin{array}{c}
\boldsymbol{v} \\ 
{A}_{1}\,\boldsymbol{u} - \xi_{0}\,\boldsymbol{v} + \kappa\,\nabla\times\boldsymbol{\varphi} - b\,\nabla\left(\theta + \alpha_{0}\,\Theta\right) \\  
\boldsymbol{\psi}  \\
{A}_{2}\boldsymbol{\varphi} - \xi\,\boldsymbol{\psi} + \kappa\,\nabla\times \boldsymbol{u} + \displaystyle\int_{0}^{+\infty}g(s)\,\Delta\boldsymbol{\eta}\ ds  \\
\Theta   \\  
k\,\Delta\theta - d\,\Theta  - b\,{\rm{div}}\,\boldsymbol{v}  \\
-\, \boldsymbol{\eta}_{s} + \boldsymbol{\psi}
\end{array} 
\right),
\end{equation*}
\begin{eqnarray*}
A_{1}\boldsymbol{u} & = & \mu\,\Delta\,\boldsymbol{u} + (\lambda + \mu + \kappa)\,\nabla\,{\rm{div}}\,\boldsymbol{u} - \kappa\nabla\times\nabla\times \boldsymbol{u},  \\
A_{2}\boldsymbol{\varphi} & = & \beta_{0}\,\Delta\,\boldsymbol{\varphi} + (\alpha + \beta)\,\nabla\,{\rm {div}}\,\boldsymbol{\varphi} - 2\,\kappa\,\boldsymbol{\varphi}
\end{eqnarray*}
and
\begin{eqnarray*}
\mathcal{D}(\mathcal{A}) & = & \left\{U\in [H^{2}(\Omega)\cap H_{0}^{1}(\Omega)]^{3}\times [H_{0}^{1}(\Omega)]^{3}\times H_{0}^{1}(\Omega)]^{3}\times [H_{0}^{1}(\Omega)]^{3} \right. \\
&  & \qquad\ \left. \times\ [ H^{2}(\Omega)\cap H_{0}^{1}(\Omega)]\times H_{0}^{1}(\Omega)\times [{\cal L}_{g} ]^{3}\,\right\},
\end{eqnarray*}
\begin{equation*}
\beta_0\,\Delta \,\boldsymbol{\varphi} + (\alpha + \beta)\,\nabla\,{\rm {div}}\,\boldsymbol{\varphi} +\displaystyle\int_{0}^{+\infty}g(s)\,\Delta\boldsymbol{\eta}\ ds\in L^{2}(\Omega)\},
\end{equation*}
where
\begin{equation}
\label{303}{\cal L}_{g} = \left\{\eta\in L_{g},\ \eta_{s}\in L_{g},\ \eta (0) = 0\right\}.
\end{equation}
Consequently, our initial boundary value problem \eqref{101}-\eqref{109} is equivalent to the Cauchy problem \eqref{302}. 

\hfill

The space $\mathcal{H}$ is equipped with a natural inner product $\langle \cdot\, ,\, \cdot \rangle_{\mathcal{H}}$, for
\begin{eqnarray*} 
U = (\boldsymbol{u},\,\boldsymbol{v},\,\boldsymbol{\varphi},\,\boldsymbol{\psi},\,\theta,\,\Theta,\,\boldsymbol{\eta}),\quad U^{*} = (\boldsymbol{u}^{*},\,\boldsymbol{v}^*,\,\boldsymbol{\varphi}^{*},\,\boldsymbol{\psi}^{*},\,\theta^{*},\,\Theta^{*},\,\boldsymbol{\eta}^{*})\in \mathcal{H}, 
\end{eqnarray*}
\begin{eqnarray*}
\lefteqn{ \langle U,\, U^{*} \rangle_{\mathcal{H}} }  \\
& = & \mu\,\displaystyle\int_{\Omega}\nabla \boldsymbol{u} (\boldsymbol{x})\cdot\nabla {\boldsymbol{u}^{*}}(\boldsymbol{x})\ d\boldsymbol{x} 
+ (\lambda + \mu + \kappa)\, \displaystyle\int_{\Omega} {\rm{div}} \boldsymbol{u} (\boldsymbol{x}){\rm{div}} {\boldsymbol{u}^{*}}(\boldsymbol{x})\ d\boldsymbol{x} 
+ \displaystyle\int_{\Omega} \boldsymbol{v} (\boldsymbol{x}){\boldsymbol{v}^* }(\boldsymbol{x})\, d\boldsymbol{x}
\end{eqnarray*}
\begin{equation*}
+\ \kappa\, \displaystyle\int_{\Omega} \left(\nabla \times \boldsymbol{u} (\boldsymbol{x}) - \boldsymbol{\varphi} (\boldsymbol{x})\right)\left(\nabla \times {\boldsymbol{u}^{*}}(\boldsymbol{x}) - {\boldsymbol{\varphi}^{*}} (\boldsymbol{x})\right)\, d\boldsymbol{x} +\kappa\, \displaystyle\int_{\Omega} \boldsymbol{\varphi} (\boldsymbol{x}) 
{\boldsymbol{\varphi}^*} (\boldsymbol{x})\ d\boldsymbol{x} 
\end{equation*}
\begin{equation*}
+\ \beta_{0}\, \displaystyle\int_{\Omega} \nabla\boldsymbol{\varphi} (\boldsymbol{x})\cdot\nabla {\boldsymbol{\varphi}^{*}} (\boldsymbol{x})\ d\boldsymbol{x} + (\alpha +\beta)\,\displaystyle\int_{\Omega} {\rm{div}}\boldsymbol{\varphi} (\boldsymbol{x}){\rm{div}} {\boldsymbol{\varphi}^{*}} (\boldsymbol{x})\ d\boldsymbol{x} + \displaystyle\int_{\Omega} \boldsymbol{\psi} (\boldsymbol{x}){\boldsymbol{\psi}^*} (\boldsymbol{x})\ d\boldsymbol{x} 
\end{equation*}
\begin{equation*}
+\displaystyle\int_{\Omega}\displaystyle\int_{0}^{+\infty}g(s)\,\nabla\boldsymbol{\eta} (\boldsymbol{x} ,s) \cdot\nabla{\boldsymbol{\eta}^*}(\boldsymbol{x},s)\,ds\, d\boldsymbol{x}  
+\, \frac{1}{\alpha_{0}}\,\displaystyle\int_{\Omega} \left(\theta (\boldsymbol{x})+ \alpha_{0}\,\Theta(\boldsymbol{x})\right) 
\left({{\theta}^{*}} (\boldsymbol{x}) + \alpha_{0}\,{{\Theta}^{*} }(\boldsymbol{x})\right)
\, d\boldsymbol{x} 
\end{equation*}
\begin{equation*}
+\  \frac{(\alpha_{0}\,d - 1)}{\alpha_{0}}\,\displaystyle\int_{\Omega} \boldsymbol{\theta} (\boldsymbol{x}) 
{\boldsymbol{\theta}^{*}} (\boldsymbol{x})\ d\boldsymbol{x} + \alpha_{0}\,k\,\displaystyle\int_{\Omega} \nabla\boldsymbol{\theta} (\boldsymbol{x}) \cdot\nabla{\boldsymbol{\theta}^{*}}(\boldsymbol{x})\ d\boldsymbol{x}.
\end{equation*}
Therefore, under the (mathematically) natural hypothesis \eqref{hyp}, $\mathcal{H}$ is a Hilbert space. Now, we claim the following well-posedness results of \eqref{302}:

\begin{theorem}
\label{existence}
Suppose that {\bf (G1)} holds true. Then, for any $U_{0} \in \mathcal{H}$, there exists a unique solution $U$ to problem \eqref{302} satisfying 
\begin{equation*}
U \in C\left(\mathbb{R}^{+};\,\mathcal{H}\right).
\end{equation*}
If moreover $U_{0} \in  \mathcal{D} (\mathcal{A})$, then
\begin{equation*}
U\in C\left(\mathbb{R}^{+};\,\mathcal{D}(\mathcal{A})\right)\cap C^{1}\left(\mathbb{R}^{+};\,\mathcal{H}\right).
\end{equation*}
\end{theorem}
\noindent 
Proof.
The proof relies on the Lumer-Philips theorem by proving that the operator $\mathcal{A}$ is dissipative and 
$Id-\mathcal{A}$ is surjective ($I\,d$ denotes the identity operator); that is $-\mathcal{A}$ is maximal monotone. So $\mathcal{A}$ is the infinitesimal generator of a $C_0$ semigroup of contraction on $\mathcal{H}$ and its domain $\mathcal{D}(\mathcal{A})$ is dense in $\mathcal{H}$. The conclusion then follows immediately (see \cite{pazy}).
\\ 
\\
By direct computations, we see that
\begin{eqnarray*}
\langle \mathcal{A}U,\, U\rangle_{\mathcal{H}} & = & -\ \xi_{0}\,\int_{\Omega}|\boldsymbol{v}|^{2}\ d\boldsymbol{x} - \xi\,\int_{\Omega}|\boldsymbol{\psi}|^{2}\ d\boldsymbol{x} - k\int_{\Omega}|\nabla\theta|^{2}\ d\boldsymbol{x} 
\end{eqnarray*}
\begin{equation}
\label{2130*}
-\, (\alpha_{0}\,d - 1) \int_{\Omega}|\Theta|^{2}\ d\boldsymbol{x} + \frac{1}{2}\int_{\Omega}\int_{0}^{+\infty}g'(s)\,|\nabla\boldsymbol{\eta}|^{2}\ ds\ d\boldsymbol{x}.
\end{equation}
Thus the operator $\mathcal{A}$ is dissipative thanks to \eqref{hyp} and the nonincreasingness of $g$. Notice that, according to \eqref{201*} and because $\boldsymbol{\eta}\in [L_{g}]^{3}$, the last term in \eqref{2130*} is well defined; indeed
\begin{equation*}
\left\vert\int_{\Omega}\int_{0}^{+\infty}g'(s)\,|\nabla\boldsymbol{\eta}|^{2}\ ds\ d\boldsymbol{x}\right\vert = - \int_{\Omega}\int_{0}^{+\infty}g'(s)\,|\nabla\boldsymbol{\eta}|^{2}\ ds\ d\boldsymbol{x} 
\end{equation*}
\begin{equation*}
\leq d_{0} \int_{\Omega}\int_{0}^{+\infty}g(s)\,|\nabla\boldsymbol{\eta}|^{2}\ ds\ d\boldsymbol{x} = d_{0}\,\Vert\boldsymbol{\eta}\Vert_{[L_{g}]^{3}}^{2} < + \infty.
\end{equation*}
\\ 
Now, we show that $I\,d - \mathcal{A}$ is surjective. Let $F = (\boldsymbol{f}_{1},\,\boldsymbol{f}_{2},\,\boldsymbol{f}_{3},\,\boldsymbol{f}_{4},\,f_{5},\,f_{6},\,\boldsymbol{f}_{7})\in \mathcal{H}$. We prove that there exists $U\in \mathcal{D}(\mathcal{A})$ satisfying
\begin{equation}
U - \mathcal{A}U=F.  \label{ZF}
\end{equation}
First, the first, third and fifth equations in \eqref{ZF} are equivalent to
\begin{equation}
\boldsymbol{v} =\boldsymbol{u} -\boldsymbol{f}_{1} ,\quad \boldsymbol{\psi} =\boldsymbol{\varphi} - \boldsymbol{f}_{3} \quad\hbox{and}\quad \Theta =\theta - f_{5}. \label{z1f1}
\end{equation}
Second, from \eqref{z1f1}, we see that the last equation in \eqref{ZF} is reduced to 
\begin{equation}
\boldsymbol{\eta}_{s} + \boldsymbol{\eta} = \boldsymbol{\varphi} + \boldsymbol{f}_{7} - \boldsymbol{f}_{3}. \label{z7*}
\end{equation}
Integrating with respect to $s$ and noting that $\boldsymbol{\eta}$ should satisfy $\boldsymbol{\eta} (0) = 0$, we get
\begin{equation}
\boldsymbol{\eta} (s) = (1 - e^{-\,s})(\boldsymbol{\varphi} -\boldsymbol{f}_{3}) + \int_{0}^{s} e^{\tau - s}\,\boldsymbol{f}_{7} (\tau)\ d\tau. \label{z7}
\end{equation}
Third, using \eqref{z1f1} and \eqref{z7}, we find that the second, fourth and sixth equations in $\eqref{ZF}$
are reduced to
\begin{equation}
\left\{
\begin{array}{ll}
\label{z5f5}
(1 + \xi_{0})\,\boldsymbol{u} - {A}_{1}\boldsymbol{u} - \kappa\,\nabla\times\boldsymbol{\varphi} + b\,(1 + \alpha_{0})\,\nabla\theta = b\alpha_{0}\,\nabla\,f_{5} + (1 + \xi_{0})\,\boldsymbol{f}_{1} +\boldsymbol{f}_{2} , \vspace{0.2cm}  \\
(1 + \xi )\,\boldsymbol{\varphi} - {A}_{2}\boldsymbol{\varphi} - g_{1}\, \Delta\boldsymbol{\varphi} - \kappa\,\nabla\times \boldsymbol{u} \\
= -\displaystyle\int_{0}^{+\infty}g(s)\displaystyle\int_{0}^{s}e^{\tau - s}\,\Delta(\boldsymbol{f}_{7}(\tau) + \boldsymbol{f}_{3})\ d\tau\ ds + (1 + \xi)\,\boldsymbol{f}_{3} + \boldsymbol{f}_{4}, \vspace{0.2cm}  \\
(1 + d)\,\theta - \kappa\,\Delta\theta + b\,{\rm{div}}\boldsymbol{u} =  (1 + d)\,f_{5}  + b\,{\rm{div}}\,\boldsymbol{f}_{1} + f_{6},
\end{array}
\right. 
\end{equation}
where 
\begin{equation*}
g_{1} = \int_{0}^{+\infty} (1 - e^{-\,s})\,g(s)\ ds.
\end{equation*}
We see that, if $\eqref{z5f5}$ admits a solution satisfying the required regularity in $\mathcal{D}(\mathcal{A})$, then \eqref{z1f1} implies that $\boldsymbol{v},$ $\boldsymbol{\psi}$ and $\Theta$ exist and satisfy the required regularity in $\mathcal{D}(\mathcal{A})$. On the other hand, \eqref{z7} implies that $\boldsymbol{\eta}$ exists and satisfies $\boldsymbol{\eta}_{s} ,\,\boldsymbol{\eta}\in [L_{g}]^{3}$; indeed, from $\eqref{z7*}$, we remark that it is enough to prove that $\boldsymbol{\eta}\in [L_{g}]^{3}$. We have 
\begin{equation*}
s\mapsto (1 - e^{-\,s})(\boldsymbol{\varphi} - \boldsymbol{f}_{3})\in [L_{g}]^{3} 
\end{equation*}
because $\boldsymbol{\varphi},\,\boldsymbol{f}_{3}\in [H_{0}^{1} (\Omega)]^{3}$. On the other hand, using the Fubini theorem and H\"older inequalities, we get
\begin{equation*}
\int_{\Omega}\int_{0}^{+\infty}g(s)\left\vert\nabla \int_{0}^{s}e^{\tau - s}\,\boldsymbol{f}_{7}(\tau)\ d\tau \right\vert^{2}ds\ d\boldsymbol{x}
\end{equation*}
\begin{eqnarray*}
& \leq & \int_{\Omega}\int_{0}^{+\infty}e^{-\,2\,s}\,g (s)\left(\int_{0}^{s}e^{\tau}\ d\tau\right)\int_{0}^{s}e^{\tau}\,\vert\nabla\boldsymbol{f}_{7}(\tau)\vert^{2}\ d\tau\ ds\ d\boldsymbol{x}   \\
& \leq & \int_{\Omega}\int_{0}^{+\infty}e^{-\,s}\,(1 - e^{-\,s})\,g (s)\int_{0}^{s}e^{\tau}\,\vert\nabla \boldsymbol{f}_{7}(\tau)\vert^{2}\ d\tau\ ds\ d\boldsymbol{x}  \\
& \leq & \int_{\Omega}\int_{0}^{+\infty}e^{-\,s}\,g(s)\int_{0}^{s} e^{\tau}\,\vert\nabla\boldsymbol{f}_{7}(\tau)\vert^{2}\ d\tau\ ds\ d\boldsymbol{x} \\
& \leq & \int_{\Omega}\int_{0}^{+\infty}e^{\tau}\,\vert\nabla \boldsymbol{f}_{7}(\tau)\vert^{2}\int_{\tau}^{+\infty}e^{-\,s}\,g (s) \,ds\ d\tau\ d\boldsymbol{x} \\
& \leq & \int_{\Omega}\int_{0}^{+\infty}e^{\tau}\,g (\tau)\,\vert\nabla\boldsymbol{f}_{7}(\tau)\vert^{2}\int_{\tau}^{+\infty}e^{-\,s}\,ds\ d\tau\ d\boldsymbol{x} \\
& \leq & \int_{\Omega}\int_{0}^{+\infty}g(\tau)\vert\nabla \boldsymbol{f}_{7}(\tau)\vert^{2}\ d\tau\ d\boldsymbol{x} \\
& \leq & \Vert \boldsymbol{f}_{7}\Vert_{[L_{g}]^{3}}^{2} < +\infty,
\end{eqnarray*}  
then 
\begin{equation*}
s\mapsto \int_{0}^{s}e^{\tau - s}\,\boldsymbol{f}_{7}(\tau)\ d\tau \in [L_{g}]^{3},
\end{equation*}
and therefore $\boldsymbol{\eta}\in [L_{g}]^{3}$. Finally, to prove that $\eqref{z5f5}$ admits a solution satisfying the required regularity in $\mathcal{D}(\mathcal{A})$, we consider the variational formulation of $\eqref{z5f5}$ and using the Lax-Milgram theorem and classical elliptic regularity arguments. This proves that $\eqref{ZF}$ has a unique solution $U\in \mathcal{D}\left( \mathcal{A}\right)$. By the resolvent identity, we have $\lambda\,Id - \mathcal{A}$ is surjective, for any $\lambda >0$ (see \cite{liu00}). Consequently, the Lumer-Phillips theorem implies that $\mathcal{A}$ is the infinitesimal generator of a linear $C_{0}$ semigroup of contractions on $\mathcal{H}$. 

\section{Asymptotic behaviour}\label{sec:stability}

In this section, we prove the main result of this paper, that is to say the stability of the micropolar thermoelastic system \eqref{302}. We start by defining the total energy and proving its nonincreasingness. Let $U_{0} \in \mathcal{H}$ and $U$ be the solution to \eqref{302}. Then, the total energy of $U$ is given by  
\begin{equation}
\label{212}E(t) = E_{w}(t) + E_{k}(t) + E_{\theta}(t),
\end{equation}
where $E_{w}$, $E_{k}$ and $E_{\theta}$ are, respectively, the potential energy, the kinetic energy and the heat conduction defined by 
\begin{align}
\lefteqn{E_{w}(t) = \frac{1}{2}\int_{\Omega}\left( \mu \,|\nabla \boldsymbol{u}|^{2}
+ (\lambda + \mu + \kappa)\,|{\rm{div}}\,\boldsymbol{u}|^{2} + \kappa \, |\nabla \times \boldsymbol{u} - \boldsymbol{\varphi}|^{2}\right)d\boldsymbol{x} }  \nonumber \\
\label{214}& +\, \frac{1}{2}\int_{\Omega}\left( \kappa \,|\boldsymbol{\varphi}|^{2}  + \beta_{0}\,|\nabla \boldsymbol{\varphi}|^{2}
+ (\alpha + \beta)\,|{\rm{div}}\,\boldsymbol{\varphi}|^{2}\right)d\boldsymbol{x} + \frac{1}{2} \int_{\Omega}\int_{0}^{+\infty}g(s)\,|\nabla\boldsymbol{\eta}|^{2}\ ds\ d\boldsymbol{x},
\end{align}
\begin{equation}
\label{215}E_{k}(t) = \frac{1}{2}\int_{\Omega}\left( |\boldsymbol{u}_{t}|^{2} +
 |\boldsymbol{\varphi}_{t}|^{2} \right)d\boldsymbol{x}
\end{equation}
and
\begin{equation}
\label{216}E_{\theta} (t) = \frac{1}{2}\int_{\Omega}\left(\frac{1}{\alpha_{0}}\,|\theta + \alpha_{0}\,\theta_{t}|^{2} + \frac{(\alpha_{0}\,d - 1)}{\alpha_{0}}\,|\theta|^{2} + \alpha_{0}\,k\,|\nabla \theta|^{2}\right)d\boldsymbol{x}.
\end{equation}
\begin{lemma}
Suppose that {\bf (G1)} holds true. Then the total energy $E$ satisfies (we use here $'$ to denote the derivative with respect to $t$)
\begin{equation*}
E' (t) = -\, \xi_{0}\,\int_{\Omega}|\boldsymbol{u}_{t}|^{2}\ d\boldsymbol{x} - \xi\,\int_{\Omega}|\boldsymbol{\varphi}_{t}|^{2}\ d\boldsymbol{x} - k\int_{\Omega}|\nabla\theta|^{2}\ d\boldsymbol{x}
\end{equation*}
\begin{equation}\label{213}
-\ (\alpha_{0}\,d - 1) \int_{\Omega}|\theta_{t}|^{2}\ d\boldsymbol{x} + \frac{1}{2}\int_{\Omega}\int_{0}^{+\infty}g'(s)\,|\nabla\boldsymbol{\eta}|^{2}\ ds\ d\boldsymbol{x} \leq 0.
\end{equation}
\end{lemma}
\noindent
Proof 
Notice that $E(t)=\frac{1}{2}\Vert U(t)\Vert_{\mathcal{H}}^{2}$. Then, according to \eqref{302},  
\begin{equation*}
E' (t)=\langle U_t (t)\, ,\, U(t)\rangle_{\mathcal{H}}=\langle \mathcal{A}U(t)\, ,\, U(t)\rangle_{\mathcal{H}} .
\end{equation*}
So, using \eqref{2130*}, we find the equality in \eqref{213}. Thanks to \eqref{hyp} and the fact that $g$ is nonincreasing, we get $E'(t)\leq 0$.  
\\
To state our stability result, we consider the following additional hypothesis on the relaxation function $g$: \\
\\
{\bf (G2)} We assume that $g(0) >0$ and there exist $d_1 >0$ and an increasing strictly convex function 
$G: \mathbb{R}^{+}\rightarrow  \mathbb{R}^{+}$ of class $C^{1}(\mathbb{R}^{+})\cap C^{2}(0,\,+\infty)$ satisfying 
\begin{equation}
\label{202}G(0) = G'(0) = 0\quad \mbox{and}\quad \lim_{t\rightarrow+\infty}G'(t) = +\infty
\end{equation}
such that
\begin{equation}
\label{203+}g'(s) \leq -\, d_{1}\,g(s),\quad\forall\;s\in \mathbb{R}^{+}
\end{equation}
or
\begin{equation}
\label{203}\int_{0}^{+\infty}\frac{g(s)}{G^{-1}\left(-\,g'(s)\right)}\ ds +\sup_{s\in\mathbb{R}^{+}}\frac{g(s)}{G^{-1}\left(-\,g'(s)\right)} < +\infty.
\end{equation}
\begin{theorem}
\label{theo:main} Assume {\bf (G1)} and {\bf (G2)} hold true. Let $U_{0} \in \mathcal{H}$ such that 
\begin{equation}
\label{203*}\eqref{203+}\,\hbox{holds}\quad\hbox{or}\quad\sup_{t\in \mathbb{R}^{+}} \int_{\Omega}\int_{t}^{+\infty}\frac{g(s)}{G^{-1}\left(-\,g'(s)\right)}\vert\,\nabla\boldsymbol{\varphi}_{0} (\boldsymbol{x},\,s - t)\vert^{2}\ ds\ d\boldsymbol{x} < +\infty.
\end{equation}
Then there exist two positive constants $c_{1}$ and $c_{2}$ such that the solution $U$ of \eqref{302} satisfies
\begin{equation}
\label{203**}E(t) \leq c_{2}\,G_{1}^{-\,1}(c_{1}\,t),\quad \forall\;t\in \mathbb{R}^{+}, 
\end{equation}
where
\begin{equation}
\label{203***} G_{1}(s) = \int_{s}^{1}\frac{1}{G_{0}(\tau)}\ d\tau\quad\mbox{and}\quad G_{0}(s) =
\begin{cases}
s & \quad\hbox{if \quad\eqref{203+}\quad \mbox{holds}}, \\
s\,G'(s) & \quad\hbox{if \quad\eqref{203}\quad \mbox{holds}}.
\end{cases}
\end{equation} 
\end{theorem}

\begin{remark}
The hypothesis \eqref{203+} implies that $g$ converges exponentially to zero at infinity. In this case, \eqref{203**} leads to the exponential stability
\begin{equation}
\label{expon} E(t) \leq c_{2}\,e^{-\,c_{1}\,t},\quad \forall\;t\in \mathbb{R}^{+} .
\end{equation} 
However, the hypothesis \eqref{203}, which was introduced by the first author in \cite{gues1}, allows $s\mapsto g(s)$ to have a decay rate at infinity arbitrarily closed to $\frac{1}{s}$. Indeed, for example, for $g(s) = q_{0}\,(1 + s)^{-\,q}$ with $q_0>0$ and $q>1$, hypothesis \eqref{203} is satisfied with $G(s)=s^{r}$, for all $r>\frac{q + 1}{q - 1}$. And then \eqref{203**} implies that 
\begin{equation*}
E(t) \leq c_{2}\,(t + 1)^{\frac{-\,1}{r - 1}},\quad \forall\;t\in \mathbb{R}^{+} .
\end{equation*}
\end{remark}

\paragraph{Proof of Theorem \ref{theo:main}.} In order to prove Theorem \ref{theo:main}, we will need to construct a Lyapunov functional equivalent to the energy $E$. For this, we will prove several lemmas with the purpose of creating negative counterparts of the terms that appears in the energy. To simplify the computations, we denote by $C$ a positive constant depending contunuously on $E(0)$ and which can be different from line to line.  
We define
\begin{equation}
\label{401}{\cal F}(t) = \int_{\Omega}\boldsymbol{u}\cdot \boldsymbol{u}_{t}\ d\boldsymbol{x} + \int_{\Omega}\boldsymbol{\varphi}\cdot \boldsymbol{\varphi}_{t}\ d\boldsymbol{x}
\end{equation}
and 
\begin{equation}
\label{402}{\cal L}(t) = E(t) + \varepsilon\,{\cal F}(t),
\end{equation}
where $\varepsilon>0$ is a small parameter to be chosen later. Next, we prove two preliminary lemmas.
\begin{lemma}
\label{lem:01}
The time derivative of the functional ${\cal F}$ defined by \eqref{401} satisfies	
\begin{equation}
\label{403}{\cal F}'(t) \leq - E(t) + C\,(E_{\theta}(t) + \,E_{k}(t)) +C\int_{\Omega}\int_{0}^{+\infty}g(s)\,|\nabla\boldsymbol{\eta}|^{2}\, ds\, d\boldsymbol{x}.
\end{equation}
\end{lemma}
Proof.  Taking the inner product in $[L^{2}(\Omega)]^{3}$ of \eqref{104} with $\boldsymbol{u}$, and the one of \eqref{207} with $\boldsymbol{\varphi}$, and adding up, we obtain
\begin{align}
{\cal F}'(t) =\ & 2\,E_{k}(t) - 2\,E_{w}(t) + b\int_{\Omega} (\theta + \alpha_{0}\,\theta_{t})\,{\rm div}\,\boldsymbol{u}\ d \boldsymbol{x} - \int_{\Omega} (\xi_{0}\,\boldsymbol{u}\,\boldsymbol{u}_{t} + \xi\, \boldsymbol{\varphi}\,\boldsymbol{\varphi}_{t})\ d \boldsymbol{x}\nonumber  \\
&\ -\, \int_{\Omega}\nabla\boldsymbol{\varphi}\int_{0}^{+\infty}g(s)\,\nabla\boldsymbol{\eta}\ ds\ d\boldsymbol{x} + \int_{\Omega}\int_{0}^{+\infty}g(s)\,|\nabla\boldsymbol{\eta}|^{2}\ ds\ d\boldsymbol{x}
 \nonumber  \\
=\ & 4\,E_{k}(t) - 2\,E(t) + 2\,E_{\theta}(t) + b\int_{\Omega}(\theta + \alpha_{0}\,\theta_{t})\,{\rm div}\,\boldsymbol{u}\ d\boldsymbol{x} - \int_{\Omega} (\xi_{0}\,\boldsymbol{u}\,\boldsymbol{u}_{t} + \xi\, \boldsymbol{\varphi}\,\boldsymbol{\varphi}_{t})\ d \boldsymbol{x} \nonumber  \\
\label{404}&\ -\int_{\Omega}\nabla\boldsymbol{\varphi}\cdot\int_{0}^{+\infty}g(s)\,\nabla\boldsymbol{\eta}\ ds\ d\boldsymbol{x} + \int_{\Omega}\int_{0}^{+\infty}g(s)\,|\nabla\boldsymbol{\eta}|^{2}\ ds\ d\boldsymbol{x}.
\end{align}
Applying the Young's and H\"older inequalities to the three two terms on the right hand side, we have, for any $\delta_{1} >0$,
there exists $C_{\delta_{1}} >0$ (depending on $\delta_{1}$) such that
\begin{equation*}
- \int_{\Omega} (\xi_{0}\,\boldsymbol{u}\,\boldsymbol{u}_{t} + \xi\, \boldsymbol{\varphi}\,\boldsymbol{\varphi}_{t})\ d\boldsymbol{x} \leq 
\delta_{1}\,E(t) + C_{\delta_{1}}\,E_{k} (t),
\end{equation*}
\begin{equation*}
b\int_{\Omega}(\theta + \alpha_{0}\,\theta_{t})\,{\rm {div}}\, \boldsymbol{u}\, d\boldsymbol{x} \leq \delta_{1} \,E(t) + C_{\delta_{1}}\,E_{\theta}(t)
\end{equation*}
 and
\begin{eqnarray}
&  & \int_{\Omega}\nabla\boldsymbol{\varphi}\int_{0}^{+\infty}g(s)\,\nabla\boldsymbol{\eta}\ ds\ d\boldsymbol{x} + \int_{\Omega}\int_{0}^{+\infty}g(s)\,|\nabla\boldsymbol{\eta}|^{2}\ ds\ d\boldsymbol{x} \nonumber \\
& \leq &\delta_{1}\,E(t) + C_{\delta_{1}}\int_{\Omega}\int_{0}^{+\infty}g(s)\,|\nabla\boldsymbol{\eta}|^{2}\ ds\ d\boldsymbol{x}.
\end{eqnarray}
Choose $\delta_{1} = \frac{1}{3}$ and combining the above three inequalities with \eqref{404}, the lemma \ref{lem:01} follows.
\begin{lemma}
\label{lem:02}
There exists a constants $\mu_{0} > 0$ such that
\begin{equation}
\label{406}-\, \mu_{0}\,E \, \leq \mathcal{F} \leq \mu_{0} \,E.
\end{equation}
\end{lemma}
\noindent 
Proof. 
From the Young's and Poincar\'e's inequalities, we have
\begin{equation*} 
\left|\mathcal{F} (t)\right| \leq \left|\int_{\Omega}\boldsymbol{u} \cdot \boldsymbol{u}_{t}\ d\boldsymbol{x}\right|  + \left|\int_{\Omega}\boldsymbol{\varphi} \cdot \boldsymbol{\varphi}_{t}\,d\boldsymbol{x}\right| \leq \mu_{0}\,E(t),
\end{equation*}
for some $\mu_{0} >0$. The conclusion of Lemma \ref{lem:02} follows straightforward .
\\
To estimate the last term of $E_{w}$, we adapt to our system a lemma introduced by the first author in \cite{gues1} and improved in \cite{gues2}. We give the proof of this lemma for the convenience of readers. 
\\
\begin{lemma}
\label{lem:020} (Lemma 3.6 \cite{gues2}) 
There exists a positive constant $c_{0}$ such that, for any $\varepsilon_{0}>0$, the following inequality holds:
\begin{equation}
\frac{G_{0}\,(\varepsilon_{0}E(t))}{\varepsilon_{0}E(t)}\displaystyle\int_{\Omega}\displaystyle\int_{0}^{+\infty }g(s)\,|\nabla\boldsymbol{\eta}|^{2}\ ds\ d\boldsymbol{x}\leq -\,c_{0}\,E^{\prime}(t)+ c_{0}\,G_{0}\,(\varepsilon_{0}E(t)).\label{3.23}
\end{equation}
\end{lemma}
\noindent 
Proof.
If \eqref{203+} holds, then we have from \eqref{213} 
\begin{equation}
\label{408}\displaystyle\int_{\Omega}\int_{0}^{+\infty}g(s)\,|\nabla\boldsymbol{\eta}|^{2}\ ds\ d\boldsymbol{x} \leq -\ \frac{1}{d_{1}}\displaystyle\int_{\Omega}\int_{0}^{+\infty}g'(t)\,|\nabla\boldsymbol{\eta}|^{2}\ ds\ d\boldsymbol{x} \leq -\ \frac{2}{d_{1}}\,E' (t).
\end{equation} 
So \eqref{3.23} holds with $c_{0} = \frac{2}{d_{1}}$ and $G_{0}(s)=s$.
\\
\\
When \eqref{203} is satisfied, we note first that, if $E (t_{0})=0$, for some $t_{0} \geq 0$, then $E(t)=0$, for all $t\geq t_{0}$, since $E$ is nonnegative and nonincreasing, and consequently, \eqref{203**} is satisfied. Thus, without loss of generality, we can assume that $E>0$ on $\mathbb{R}^+$.
\\
\\
Because $E$ is nonincreasing, we have
\begin{eqnarray*}
\int_{\Omega}\,|\nabla\boldsymbol{\eta}|^{2}\ d\boldsymbol{x}
& \leq & 2\left(\int_{\Omega}\,|\nabla\boldsymbol{\varphi}(\boldsymbol{x},\,t)|^{2}\ d\boldsymbol{x} + \int_{\Omega}\,|\nabla\boldsymbol{\varphi}(\boldsymbol{x},\,t - s)|^{2}\ d\boldsymbol{x} \right)  \\
& \leq & C\,E(0) + 2\int_{\Omega}|\nabla\boldsymbol{\varphi}(\boldsymbol{x},\,t - s)|^{2}\ d\boldsymbol{x}  \\ 
& \leq & \left\{
\begin{array}{ll}
C\,E(0) & \quad\hbox{if} \quad 0\leq s\leq t, \\
C\,E(0) + 2\int_{\Omega}|\nabla\boldsymbol{\varphi}_{0}\, (\boldsymbol{x},\,s - t)|^{2}\ d\boldsymbol{x}&\quad\hbox{if}\quad s>t\geq 0
\end{array} :=M(t,\,s), 
\right.
\end{eqnarray*}
so we conclude that 
\begin{equation}
\int_{\Omega}\,|\nabla\boldsymbol{\eta}|^{2}\ d\boldsymbol{x} \leq M(t,\,s),\quad \forall \;t,\,s\in \mathbb{R}^{+}. \label{3.26} 
\end{equation} 
Let $\tau_{1}(t,\,s),\ \tau_{2}(t,\,s)>0$ (which will be fixed later on), $\varepsilon_{0} >0$ and $K(s)={\frac{s}{{G^{-\,1}(s)}}}$, for $s>0$, and $K(0)=0$, since {\bf (G2)} implies that  
\begin{equation*}
\lim_{s\to 0^{+}}\dfrac{s}{G^{-1}(s)}=\lim_{\tau\to 0^{+}}\dfrac{G(\tau)}{\tau}=G'(0)=0.
\end{equation*}
The function $K$ is nondecreasing. Indeed, the fact that $G^{-1}$ is concave and $G^{-1} (0)=0$ implies that, for any $0\leq s_{1} < s_{2}$,
\begin{eqnarray*}
K(s_{1}) 
&=& \frac{s_1}{G^{-\,1}\left(\frac{s_{1}}{s_{2}}s_{2} + \left(1 - \frac{s_{1}}{s_{2}}\right)0\right)}  
\\
& \leq & 
\frac{s_1}{ 
\frac{s_1}{s_2}\,G^{-1}(s_{2}) 
+ \left( 1 - \frac{s_1}{s_2} \right)
G^{-1}(0)
} = 
\frac{ s_2 }{ G^{-1} (s_2)} = K(s_2 ).
\end{eqnarray*} 
Then, using \eqref{3.26},
\begin{equation*}
K\left(-\,\tau_{2}(t,\,s)\,g^{\prime}(s)\int_{\Omega}|\nabla\boldsymbol{\eta}|^{2}\ d\boldsymbol{x}\right) \leq K\left(-\,M(t,\,s) \,\tau_{2}(t,\,s)\,g^{\prime}(s)\right),\quad\forall\;s\in \mathbb{R}^{+}. 
\end{equation*}
Using this inequality, we arrive at 
\begin{eqnarray*}
\lefteqn{\int_{0}^{+\infty}g(s)\int_{\Omega}|\nabla\boldsymbol{\eta}|^{2}\ d\boldsymbol{x} ds  }  \\
& = & \dfrac{1}{G'(\varepsilon_{0}\,E(t))}\int_{0}^{+\infty}\dfrac{1}{\tau_{1}(t,\,s)}G^{-\,1}\left(-\,\tau_{2}(t,\,s)\,g'(s)\int_{\Omega}|\nabla\boldsymbol{\eta}|^{2}\ d\boldsymbol{x}\right) \\
&  & \mbox{} \times\dfrac{\tau_{1}(t,\,s)\,G'(\varepsilon_{0}\,E(t))\,g(s)}{-\,\tau_{2}(t,\,s)\,g'(s)}\ K\left(-\,\tau_{2}(t,\,s)\,g'(s)\int_{\Omega}|\nabla\boldsymbol{\eta}|^{2}\ d\boldsymbol{x}\right)ds   \\
& \leq & \dfrac{1}{G'(\varepsilon_{0}\,E(t))}\int_{0}^{+\infty}\dfrac{1}{\tau_{1}(t,\,s)}\,G^{-\,1}\left(-\,\tau_{2}(t,\,s)\,g'(s)\int_{\Omega}|\nabla\boldsymbol{\eta}|^{2}\ d\boldsymbol{x}\right)  \\
&  & \mbox{} \times\dfrac{\tau_{1}(t,\,s)G'(\varepsilon_{0}\,E(t))\,g(s)}{-\,\tau_{2}(t,\,s)\,g'(s)}\ K\left(-\,M(t,\,s)\,\tau_{2}(t,\,s)\,g'(s)\right)ds   \\
& \leq & \dfrac{1}{G' (\varepsilon_{0}\,E(t))}\int_{0}^{+\infty}\dfrac{1}{\tau_{1}\,(t,\,s)}\,G^{-1}\left(-\,\tau_{2}(t,\,s)\,g'(s)\int_{\Omega}|\nabla\boldsymbol{\eta}|^{2}\ d\boldsymbol{x}\right)  \\
&  & \mbox{} \times\dfrac{M(t,\,s)\,\tau_{1}(t,\,s)\,G'(\varepsilon_{0} \,E(t))\,g(s)}{G^{-\,1}\,(-\,M(t,\,s)\,\tau_{2}(t,\,s)\,g'(s))}\ ds.
\end{eqnarray*}
Let $G^{*} (s)=\sup_{\tau\in \mathbb{R}_{+}} \{s\,\tau - G (\tau)\}$, for $s\in \mathbb{R}^{+}$, denote the dual function of $G$. Thanks to {\bf (G2)}, we see that
\begin{equation*}
G^{*}(s) = s\,(G')^{-1} (s) - G ((G')^{-1} (s)),\quad\forall\;s\in \mathbb{R}^{+} .
\end{equation*}
Using Young's inequality: $s_{1}\,s_{2}\leq G(s_{1}) + G^{*}(s_{2})$, for 
\begin{equation*}
s_{1} = G^{-1}\left(-\,\tau_{2}(t,\,s)\,g^{\prime}(s)\int_{\Omega}|\nabla\boldsymbol{\eta}|^{2}\ d\boldsymbol{x}\right)\quad\mbox{and}\quad s_{2} ={\dfrac{{M(t,\,s)\,\tau_{1}(t,\,s)\,G^{\prime}(\varepsilon_{0}\,E(t))\,g(s)}}{{G^{-\,1}(-\,M(t,\,s)\,\tau_{2}(t,\,s) g^{\prime }(s))}}}, 
\end{equation*}
we get 
\begin{eqnarray*}
\lefteqn{\int_{0}^{+\infty}g(s)\int_{\Omega}|\nabla\boldsymbol{\eta}|^{2}\ d\boldsymbol{x}\ ds \leq  \dfrac{1}{G^{\prime}(\varepsilon_{0} \,E(t))}\int_{0}^{+\infty}\dfrac{{-\,\tau_{2}(t,\,s)}}{{\tau_{1}(t,\,s)}}g^{\prime}(s)\int_{\Omega}|\nabla\boldsymbol{\eta}|^{2}\ d\boldsymbol{x}\ ds } \\
&  & \mbox{} +\ \dfrac{1}{G^{\prime}(\varepsilon_{0}\,E(t))}\int_{0}^{+\infty}\dfrac{1}{\tau_{1}(t,\,s)}\ G^{*}\left({\frac{{M(t,\,s)\,\tau_{1}(t,\,s)\,G^{\prime}(\varepsilon_{0}\,E(t))\,g(s)}}{{G^{-1}(-\,M(t,\,s)\,\tau_{2}(t,\,s)\,g^{\prime}(s))}}}\right)ds. 
\end{eqnarray*}
Using the fact that $G^{*}(s)\leq s\,(G')^{-\,1} (s)$, we get 
\begin{eqnarray*}
\int_{0}^{+\infty}g(s)\int_{\Omega}|\nabla\boldsymbol{\eta}|^{2}\ d\boldsymbol{x}\ ds  \leq  \dfrac{-\,1}{G^{\prime }(\varepsilon_{0}\,E(t))}\int_{0}^{+\infty}\dfrac{{\tau_{2}(t,\,s)}}{{\tau_{1}(t,\,s)}}\ g^{\prime }(s)\int_{\Omega}|\nabla\boldsymbol{\eta}|^{2}\ d\boldsymbol{x}\ ds 
\end{eqnarray*}
\begin{eqnarray*}
\mbox{} + \int_{0}^{+\infty}{\frac{{M(t,\,s)\,g(s)}}{{G^{-\,1}(-\,M(t,\,s)\,\tau_{2}(t,\,s)\,g^{\prime}(s))}}}\,(G')^{-\,1}\left({\frac{{M(t,\,s)\,\tau_{1}(t,\,s)\,G^{\prime}(\varepsilon_{0}\,E(t))\,g(s)}}{{G^{-\,1}(-\,M(t,\,s)\,\tau_{2}(t,\,s)\,g^{\prime}(s))}}}\right)ds. 
\end{eqnarray*}
Then, using the fact that $(G')^{-1}$ is nondecreasing and choosing $\tau_{2}(t,\,s) = {\frac{{1}}{{M(t,\,s)}}}$, we get 
\begin{eqnarray*}
\lefteqn{\int_{0}^{+\infty}g(s)\int_{\Omega}|\nabla\boldsymbol{\eta}|^{2}\ d\boldsymbol{x}\ ds } \\
& \leq & \dfrac{-\,1}{G^{\prime}(\varepsilon_{0}\,E(t))}\int_{0}^{+\infty}\dfrac{1}{M(t,\,s)\,\tau_{1}(t,\,s)}\,g^{\prime}(s)\int_{\Omega}\,|\nabla\boldsymbol{\eta}|^{2}\ d\boldsymbol{x}\ ds  \\
&  & \mbox{} + \int_{0}^{+\infty}{\dfrac{{M(t,\,s)\,g(s)}}{{G^{-\,1}(-\,g^{\prime}(s))}}}\ (G')^{-1}\left(m_{1}\,M(t,\,s)\,\tau_{1}(t,\,s)\,G^{\prime}(\varepsilon_{0}\,E(t))\right)ds,  
\end{eqnarray*}
where $m_{1} = \sup_{s\in \mathbb{R}^{+}}{\frac{{g (s)}}{{G^{-1}(-\,g^{\prime}(s))}}} <+\infty$ ($m_{1}$ exists according to \eqref{203}).
Due to \eqref{203} and the restriction on $\boldsymbol{\varphi}_{0}$ in \eqref{203*}, we have
\begin{equation*}
\sup_{t\in \mathbb{R}^{+}}\displaystyle\int_0^{+\infty}{\frac{{M(t,\,s)\,g(s)}}{{G^{-\,1}(-\,g^{\prime}(s))}}}\ ds=:m_{2} <+\infty .
\end{equation*} 
Therefore, choosing $\tau_{1}(t,\,s)= {\frac{{1}}{{m_{1}\,M(t,\,s)}}}$ and using \eqref{213}, we obtain 
\begin{eqnarray*}
\lefteqn{\int_{0}^{+\infty}g(s)\displaystyle\int_{\Omega}|\nabla\boldsymbol{\eta}|^{2}\ d\boldsymbol{x}\ ds }  \\
& \leq & \dfrac{{-\,m_{1}}}{{G^{\prime}(\varepsilon_{0}\,E(t))}}\int_{0}^{+\infty}g^{\prime}(s)\int_{\Omega}|\nabla\boldsymbol{\eta}|^{2}\ d\boldsymbol{x}\ ds + \varepsilon_{0}\,E(t)\int_{0}^{+\infty}{\frac{{M(t,\,s)\,g(s)}}{{G^{-\,1}(-\,g^{\prime}(s))}}}\ ds  \\
& \leq & \dfrac{{-\,2\,m_{1}}}{{G^{\prime}(\varepsilon_{0}\,E(t))}}\,E^{\prime}(t) + \varepsilon_{0}\,m_{2}\,E(t),   
\end{eqnarray*}
which gives \eqref{3.23} with $c_{0} = \max \{2\,m_{1},\,m_{2}\}$ and $G_{0}(s) = s\,G' (s)$. 
\\
We are now ready to prove the main stability result \eqref{203**}. First, we observe that, from \eqref{216} and the Poincar\'e's inequality, the heat energy $E_{\theta}$ is controlled by the sum of $|\theta_{t}|^{2}$ and $|\nabla\theta|^{2}$; that is
\begin{equation*}
E_{\theta} (t) \leq C\int_{\Omega}\left(|\theta_{t}|^{2} + |\nabla\theta|^{2}\right) d\boldsymbol{x}\ ds.
\end{equation*}
Therefore
\begin{equation}
\label{407}
-\,\kappa\,\int_{\Omega}|\nabla\theta|^{2}\ d\boldsymbol{x} - (\alpha_{0}\,d - 1)\,\int_{\Omega}|\theta_{t}|^{2}\ d\boldsymbol{x} \leq -\ C\,E_{\theta}(t).
\end{equation}
Now, we differentiate $\cal{L}$ from \eqref{402} with respect to time and use \eqref{403}, together with the dissipation of energy \eqref{213} and the above estimate \eqref{407}, we obtain, for some positive constants $C_{0},\,C_{1}$ and $C_{2}$,
\begin{eqnarray*}
{\cal L}' (t) & \leq & -\, \left(\xi_{0} - \frac{\varepsilon\,C_{0}}{2}\right)\int_{\Omega}|\boldsymbol{u}_{t}|^{2}\ d\boldsymbol{x} - \left(\xi - \frac{\varepsilon\,C_{0}}{2}\right)\int_{\Omega}|\boldsymbol{\varphi}_{t}|^{2}\ d\boldsymbol{x} 
\\
&  & -\ (C_{2} - \varepsilon\,C_{0})\,E_{\theta}(t) - \varepsilon\,E(t) +\varepsilon\,C_{0}\int_{\Omega}\int_{0}^{+\infty}g(s)\,|\nabla\boldsymbol{\eta}|^{2}\ ds\ d\boldsymbol{x}.
\end{eqnarray*}
Choosing 
\begin{equation}
\label{3040}
0<\varepsilon <\min\left\{\frac{2\,\xi_{0}}{C_{0}},\ \frac{2\,\xi}{C_{0}},\ \frac{C_{2}}{C_{0}}\right\},
\end{equation}
we obtain
\begin{equation}
\label{304}
{\cal L}' (t)  \leq  -\,\varepsilon\,E(t) + C\int_{\Omega}\int_{0}^{+\infty}g(s)\,|\nabla\boldsymbol{\eta}|^{2}\ ds\ d\boldsymbol{x}.
\end{equation}
Multiplying \eqref{304} by $\frac{G_{0}(\varepsilon_{0}\,E)}{\varepsilon_{0}\,E}$ and combining \eqref{3.23}, we find
\begin{equation}
\label{304**}
\frac{G_{0}(\varepsilon_{0}\,E(t))}{\varepsilon_{0}\,E(t)}\,{\cal L}' (t)  + c_{0}\,C\,E'(t) \leq -\,\left(\frac{\varepsilon}{\varepsilon_{0}} - c_{0}\,C \right)G_{0}(\varepsilon_{0}\,E(t)).
\end{equation}
Now we define our Lyapunov functional ${\cal R}$ by  
\begin{equation}
\label{304*}
{\cal R} = \tau_{0} \left(\frac{G_{0}(\varepsilon_{0}\,E)}{\varepsilon_{0}\,E}{\cal L} + c_{0}\,C\,E\right),
\end{equation}
where $\tau_{0}$ is positive constant that will be choose later. On the other hand, thanks to Lemma \ref{lem:02}, we get
\begin{equation}
\label{410} (1 - \varepsilon\,\mu_{0})\,E(t) \leq {\cal L}(t) \leq (1 + \varepsilon\, \mu_{0} )\,E(t).
\end{equation}
Choosing $\varepsilon$ such that \eqref{3040} holds and $\varepsilon <\frac{1}{\mu_{0}}$, we see that ${\cal L}$ and $E$ are equivalent. Because $E$ is nonincreasing and $G$ is convexa, then $\frac{G_{0} (\varepsilon_{0}\,E)}{\varepsilon_{0}\,E}$ is nonincreasing, and therefore, the \eqref{304**} and \eqref{304*} lead to 
\begin{equation}\label{3041}
{\cal R}' (t)  \leq -\,\tau_{0}\,\left(\frac{\varepsilon}{\varepsilon_{0} } - c_{0}\,C \right)G_{0}(\varepsilon_{0}\,E(t)).
\end{equation}
Moreover, recalling that $\frac{G_{0}(\varepsilon_{0}\,E)}{\varepsilon_{0}\,E}$ is nonincreasing and using \eqref{410}, we obtain 
\begin{equation}
\label{4101} \tau_{0}\,c_{0}\,C\,E(t) \leq {\cal R}(t) \leq \tau_{0}\left[(1 + \varepsilon\,\mu_{0})\,\frac{G_{0}(\varepsilon_{0}\,E(0))}{\varepsilon_{0}\,E(0)} + c_{0}\,C\right]\,E(t).
\end{equation}
By choosing $0<\varepsilon_{0} <\frac{\varepsilon}{c_{0}\,C}$, we deduce from \eqref{3041} and \eqref{4101} that ${\cal R}$ is equivalent to $E$ and satisfies  
\begin{equation}
\label{3042}
{\cal R}' (t)  \leq -\,\tau_{0}\,C\,G_{0}(\varepsilon_{0}\,E(t)).
\end{equation}
Thus, for $\tau_{0} >0$ such that
\begin{equation*}
{\cal R} \leq \varepsilon_{0}\,E \quad\hbox{and}\quad {\cal R}(0) \leq 1,
\end{equation*} 
we get, for $c_{1} = \tau_{0}\,C$,
\begin{equation}
\label{3.30}
{\cal R}^{\prime} (t) \leq -\ c_{1}\,G_{0}( {\cal R}(t)).
\end{equation}
Then \eqref{3.30} implies that $(G_{1}({\cal R}))' \geq c_{1}$, where $G_{1}$ is defined in \eqref{203***}. So, a direct integrating gives 
\begin{equation}
\label{3.300}
G_{1}({\cal R} (t)) \geq c_{1}\,t + G_{1}({\cal R} (0)).
\end{equation}
Because ${\cal R} (0)\leq 1$ and $G_{1}$ is decreasing, we obtain $G_{1} ({\cal R} (t)) \geq c_{1}\,t$ which implies that ${\cal R} (t)\leq G_{1}^{-1} (c_{1}\,t)$. Finally, from the equivalence of ${\cal R}$ and $E$, the result \eqref{203**} follows and the proof of Theorem \ref{theo:main} is complete.

\section{Comments and issues}

{\bf Comment 1.} If $g=0$ (that is only the linear frictional damping $\xi\,\boldsymbol{\varphi}_{t}$ is considered on \eqref{207}), then \eqref{203**} is reduced to the exponential stability estimate \eqref{expon}.
\\ 
\\
{\bf Comment 2.} The results of this paper can be generalized to the case of nonlinear damping; that is, $\xi_{0}\,\boldsymbol{u}_{t}$ and $\xi\,\boldsymbol{\varphi}_{t}$ are replaced by $h_{1}(\boldsymbol{u}_{t})$ and $h_{2}(\boldsymbol{\varphi}_{t})$, where $h_{i}: \mathbb{R}^{3}\to \mathbb{R}$ are given functions satisfying some hypotheses; see, for example, \cite{gues3} (for coupled Timoshenko systems).
\\
\\
{\bf Comment 3.} Our results hold true when \eqref{207} is controlled only via the infinity memory (that is $\xi=0$). In this case, the second integral in \eqref{401} is replaced by 
\begin{equation}
\label{61}
-\int_{\Omega}\boldsymbol{\varphi}_{t}\int_{0}^{+\infty}g(s)\,\boldsymbol{\eta}\ ds\ d\boldsymbol{x}.
\end{equation}
For more details, see, for example, \cite{gues1}, \cite{gues3} and \cite{gues2} (for Timoshenko and abstract systems).  \\
\\
{\bf Comment 4.} It is possible to consider different relaxation functions $g_{1}$, $g_{2}$ and $g_{3}$ (instead of $g$) satisfying the hypotheses ({\bf G1}) and ({\bf G2}); that is $g\Delta\boldsymbol{\varphi}$ is replaced by $\left(g_{1}\,\Delta\varphi_{1},\ g_{2}\,\Delta\varphi_{2},\ g_{3}\,\Delta\varphi_{3}\right)$. See \cite{gubc} (for \eqref{102} with $\kappa =\xi =0$).
\\
\\
{\bf Comment 5.} The equation \eqref{104} can be controlled via an infinite memory 
\begin{equation}
\label{62}
\int_{0}^{+\infty}f(s)\,\Delta\boldsymbol{u} (\boldsymbol{x},\,t - s)\ ds
\end{equation}
instead of the linear damping $\xi_{0}\,\boldsymbol{u}_{t}$, where $f\,: \mathbb{R}^{+}\to \mathbb{R}^{+}$ is a given relaxation function satisfying the same hypotheses as $g$. To prove the well-posedness results, we introduce a second variable $\boldsymbol{z}$ similar to $\boldsymbol{\eta}$ given by 
\begin{equation}
\label{63}
\boldsymbol{z}(\boldsymbol{x},\,t,\,s) = \boldsymbol{u} (\boldsymbol{x},\,t) - \boldsymbol{u} (\boldsymbol{x},\,t - s).
\end{equation}
We define its space $[L_{f}]^{3}$ as $[L_{g}]^{3}$ and do some logical modifications. For the stability result, we replace the first integral in \eqref{401} by 
\begin{equation}\label{63}
-\int_{\Omega}\boldsymbol{u}_{t}\int_{0}^{+\infty}f(s)\,\boldsymbol{z}\ ds\ d\boldsymbol{x}.
\end{equation}  
{\bf Comment 6.} In fact, our results hold true also when $\alpha_{0}\,d - 1 =0$. In this case, using Poincar\'e's and Young's inequalities, we see that, for any $0<\varepsilon < \frac{\alpha_{0}\,C_{p}}{C_{p} + \kappa\,\alpha_{0}^{2}}$ ($C_{P}$ is the Poincar\'e's constant),
\begin{eqnarray*}
E_{\theta} (t) & \geq & \frac{1}{2}\int_{\Omega}\left(\left[\left(\frac{1}{\alpha_{0}} - \frac{1}{\varepsilon}\right)C_{P} + \kappa\,\alpha_{0}\right]\,|\nabla\theta|^{2} + \left(\alpha_{0} - \varepsilon\right)\,|\theta_{t}|^{2}\right) d\boldsymbol{x} \\
& \geq & 
C\int_{\Omega} \left(|\nabla\theta|^{2} + |\theta_{t}|^{2}\right) d\boldsymbol{x}.
\end{eqnarray*}
So $E_{\theta}$ is still definite positive quadratic form. Consequently, the well-posedness result is still valid. 
For the stability, we have only to add to the definition of ${\cal F}$ in \eqref{401} the term 
\begin{equation}
\label{64}
-\int_{\Omega}\theta\,\theta_{t}\ d\boldsymbol{x}.
\end{equation} 
When $\alpha_{0}\,d - 1 <0$, the situation is more complicate because the energy is not necessarily nonincreasing (see \eqref{213}). 
\\
\\
{\bf Comment 7.} The last comment concerns the stability of our system when only one equation from \eqref{104} and \eqref{207} is controlled (via a frictional damping or infinite memory); that is $\xi = g = 0$ or $\xi_{0} = f  = 0$. So one of the equations \eqref{104} and \eqref{207} is indirectly controlled by the other one via the coupling terms. Probably, the system is still stable but maybe with a weaker decay rate that the one given by \eqref{203**}. This question will be the subject of a future work. 
\\
\\
{\bf Acknowledgment.} 
%This work was initiated during the visit in July-August 2017 of the first author to LNCC and RJ university, Brazil, and B\'{\i}o-B\'{\i}o and Concepci\'{o}n universities, Chile. 
The first author thanks LNCC, RJ, B\'{\i}o-B\'{\i}o and Concepci\'{o}n universities for their kind support and hospitality.
The third author thanks   the support of  FONDECYT grant no. 1180868, and ANID-Chile through the project {\sc Centro de Modelamiento Matem\'atico} (AFB170001)  of the PIA Program: Concurso Apoyo a Centros Cient\'\i ficos y Tecnol\'ogicos de Excelencia con Financiamiento Basal.

\end{document}